# SEGMENTATION OF THE HOMOGENEITY OF A SIGNAL USING A PIECEWISE LINEAR RECOGNITION TOOL

Joseph Morlier

LRBB mixed unit CNRS/INRA/Bordeaux 1

69 route d'Arcachon

33612 Cestas Gazinet, France

Tel: +335517122820, Fax: +33556680713

jmorlier@lrbb.u-bordeaux.fr

**Abstract**

In this paper a new method of detection of homogeneous zones and singularity parts of a 1D signal is proposed. The entropy function is used to transform signal in piecewise linear one. The multiple regression permits to detect lines and project them in the Hough parameters space in order to easily recognise homogeneous zones and abrupt changes of the signal. Two application examples are analysed, the first is a classical fractal signal and the other is issued from a dynamic mechanical study.

**Keywords**: Piecewise affine, line detection, Hough parameters, entropy, segmentation, homogeneity.

# 1 Introduction

In this paper a method for recognising homogeneous and singular parts of a signal is proposed. The detection of piecewise linear signals (or PieceWise Affine PWA) is an interesting domain of the pattern recognition sciences. Interpolation and least squares methods have been widely used [3], and recently randomised Hough transform succeed to determine line segments from serial samples in real time [6]. The first application of the Hough transform [12] is to detect lines in a binary image using mappings between points on the Euclidean plane and parameters which define lines on the plane. More recently learning techniques like neural networks are also interesting to estimate piecewise linear regression in 2D and 3D [2]. Moreover adaptive linear regression approachs are also popular and particulary in automatic target recognition. The area where the signatures have large variance is partitioned into smaller intervals to maintain high accuracy of the estimate; in the low variance area, a large interval is taken to reduce the computation cost of the classification algorithm which uses this piecewise linear template model [9].

The aim of our work is twofolded. First taking advantage of the mathematical formulation of the entropy function to transform any signal in a piecewise linear one [1]. Secondly, to create a simple algorithm of multiple regression of piecewise linear signal which permits to extract any affine line using a least-mean-squares method on recursive sampled data. Thus the projection of these detected lines in the Hough parameters space (position, length, slope) allows to easily detect homogeneous zone and abrupt changes, in other words to segment or classify the signal in two categories. These techniques are applied in this paper on two examples, first a synthesis fractal signal in order to identify the Hurst exponent [5] and the other is a dynamic modeshape data of a damaged beam with the aim of characterising the damage [4, 7, 8].

# 2 Entropy function

Entropy function have been used in soil diagraphy signal analysis in order to segment it in homogeneous zones. The concept of the thermodynamics of curves is used to analyze change in the variability of time or spatial series [1].

Considering $X(z)$, a random signal, diffentiable and supposing the absolute value of the first derivative is stationnary (order 2).

Thus, its average is constant:
$$m|x'|(z) = E[|X'|] = m|x'|$$
(1.)

its standard deviation is constant:
$$\sigma^2|x'|(z) = E[(|X'| - m|x'|)^2] = 2|x'|$$
(2.)

It is demonstrated that for time and spatial series the entropy of a curve is the slope of the cumulative sum of absolute differences. Then the entropy function H(z) of the signal X(z) can be written by:
$$H(z) = \int_0^z |x'|(u)\ du$$
(3.)

So, H(z) is linear in average:
$$E[H(z)] = E[\int_0^z |x'|(u)\ du]$$
$$= \int_0^z E[X'(u)]du = mH(z) = z\ m|x'|$$
(4.)

Thus entropy function is able to transform any signal in piecewise linear one. More specifically, the entropy of a curve makes it possible to divide a nonstationary random field, because of a change in variance, into subdomains where data are said to be stationary. It can be shown that the presence of a linear or quadratic trend is without effect on the localization of the stationary subdomains.

The next numerical simulation aims at demonstrating the efficiency of this tool on a random Gaussian signal.

Considering $y(x) = 0.1*\text{rand}(1,100)$, a random signal of noise standard deviation equal to 0.1.

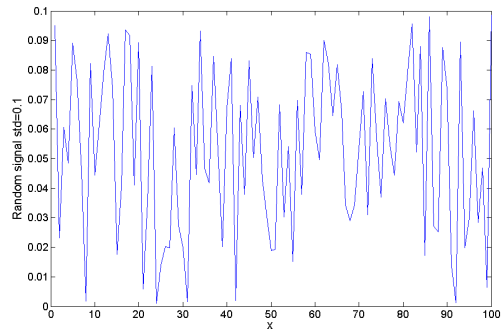

**Figure** 1. Random signal

The mean of the absolute value of the first derivative (Eq.1) is equal to the slope of the entropy function (0.0337 close to 0.0316) fo the Figure 2 below.

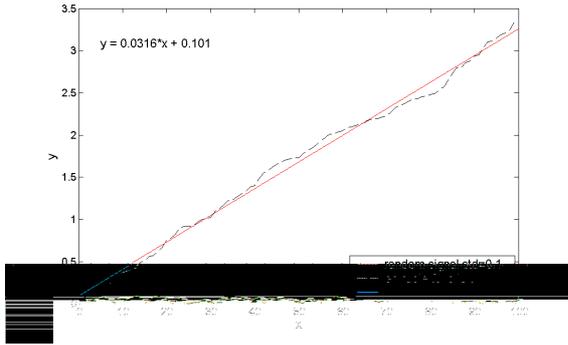

**Figure** 2. Entropy of y(x) and the smooth function of y(x)=0.0316x+0.101

## 3 Line detection and characterization

In this section, we develop a recursive multiple regression using the coefficient of correlation ($R^2$) of the sampled fit function as pertinent condition to end the line detection algorithm. Expert system should enhance the procedure to automate the choice of this coefficient value which increases as the experimental noise decreases.

The method has been tested firstly with a simulated signal constructed with a set of connecting lines segments. Considering $f(x)$ as:

$$f(x) = \begin{cases} -x & if -4 \leq x \leq 0 \\ x & if \ 0 < x < 2 \\ 2+3x & if \ 2 \leq x \leq 4 \end{cases}$$

(5.)

The estimation of the detected line are estimated with optimal $R^2=0.998$:

$$f(x) = \begin{cases} -0.99363x + 0.06774 & if -4 \leq x \leq 0 \\ 0.99176x + 0.055183 & if \ 0 < x < 2 \\ 2.0484 + 2.9992x & if \ 2 \leq x \leq 4 \end{cases}$$

(6.)

The Hough parameters are obtain using a one to one mapping between detected lines described by pairs of tangent a and the intercept b and line length L (at position p) and radius formed with the x absciss defined the slope as :

$$\alpha = \frac{180}{\pi} \tan^{-1}(a).$$

(7.)

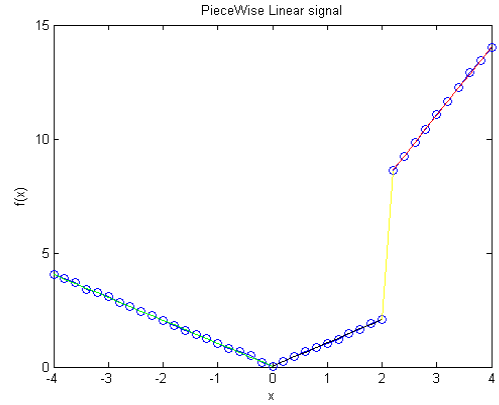

**Figure** 3. Piecewise linear signal (o) and detected lines (colored line)

The Figure 3 above represents the 4 detected lines, the coefficient estimated permits to compute the Hough parameters (Eq. 7).

The Figure 4 highlight the behaviour of the piecewise affine: 3 homogeneous zones (large length at position 0, 2 and 4 and one singularity (2 points joining affine lines 2 and 3).

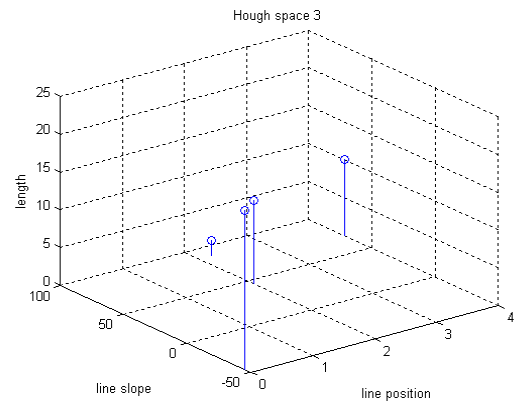

**Figure** 4. Representations of the lines (3 homogeneous parts, one singularity) projected in the Hough parameters space

The interest is poor for simple piecewise signal, but becomes more interessant when complex signals are transformed in piecewise linear using entropy and then analysed with multiple regression in order to classify the detected line in respectivly homogeneous and singularity parts.

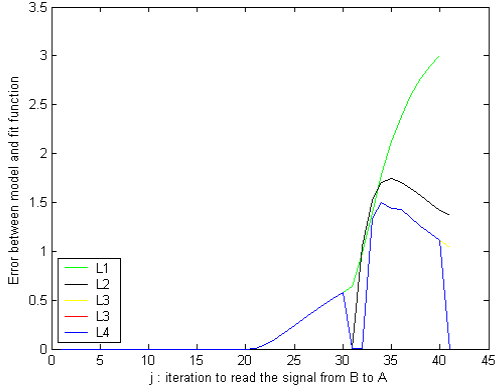

**Figure** 5. Errors of fitting for each detected lines from point B(j=41) to A(j=0) at each iteration

The results (Fig. 5) of the fit algorithm demonstrate the abilitie of the multiple regression algorithm to converge after 22 itérations for line L1, 10 for line L2, 9 to isolate the joint line L3 and 8 to complete the last line to estimate.

To study the noise robustness, we add a random gaussian noise with zero mean of different variance to the signal as:

$$f(x) = \begin{cases} -x + \varepsilon & if -4 \leq x \leq 0 \\ x + \varepsilon & if\ 0 < x < 2 \\ 2 + 3x + \varepsilon & if\ 2 \leq x \leq 4 \end{cases}$$

(8.)

Then the user choose the optimal stopping criteria that respect the 4 lines projection in the Hough space (no false lines detected).

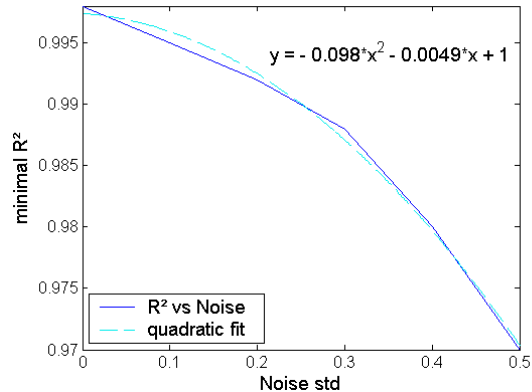

**Figure** 6. Optimal $R^2$ vs noise standard deviation added

The Figure 6 above answer us about the optimal $R^2$ behaviour to choose function of the noise standard deviation added.
Moreover it shows an estimation of the line parameters function of noise. Small variations occurs until false line were detected at noise std=0.5. We can conclude that results of the regression are still satisfying for low noise injected. It means the Hough parameters should vary in reasonable value and the interpretation should still be easy.

## 4  Segmentation algorithm

The algorithm of classification of homogeneous zones of the signal is described below:
The user should enter the estimated number of line to be detected in order to limit the number of loop in the detection scheme.
The choise of optimal $R^2$ is critical to the accuracy and efficiency of the signal classification.

```
1 :  begin procedure of signal (x,y) classification
2 :     assignate appropriate value to
        B=length(ye)
        A=1
        j=B+1
        [xe,ye]=entropy(x,y)
3 :     while R² <Rm²  and j>A+1 do
4 :     begin
             j=j-1
             yeᵢ=ye(A:j)
             xeᵢ=xe(A:j)
             [aᵢ,bᵢ]=lms(xeᵢ,yeᵢ)
             R²=corrcoef(xeᵢ,yeᵢ)²
             ylᵢ=aᵢ.xᵢ+bᵢ
             errorᵢ=mean|yeᵢ-ylᵢ|
             αᵢ = tan⁻¹(aᵢ)
             Lᵢ=length(yeᵢ)
             A=j
        end
     end
```

**Figure** 3. Algorithm of segmentation

Line 2 assign important values A and B the index of the segmented data and j the variable which reads the signal from point B to A. The entropy function is the slope of the cumulative sum of absolute differences resulting a piecewise afine signal.
The loop of multiple regressions begins at line 3, the stopping criteria is defined by the user (function of number of line to recognize and noise in the signal). For each segmented signal satisfying the conditon $R^2<Rm^2$, we use least mean squares (lms) algorithm to fit the segmented line and extract the coefficient $a_i$ and $b_i$. Then the coefficient of correlation is estimated (corrcoef) and the Hough parameters $\alpha_i$, $L_i$ and error between model and segmented data ($error_i$) are computed. The next segmented data begins with A=j.

## 5 Experimental Results

Firstly our method is validated on a Fractional Brownian motion (FBm), a popular synthesis fractal signal, which is time varying random function with stationary, Gaussian distributed, and statistically self-affined increment. Practically FBm exists for example in the fracture roughness surface of materials. The fractal dimension D describes how many pieces of a set resolved as the resolution scales decreases. It informs us about the signal complexity. The relationship between measurement (of the pieces) and measuring scales r has a power law property with the following curve.

$$C(r) = \alpha r^D$$
(9.)

Here $\alpha$ is constant and $C(r)$ is the number of pieces of a size r covering the signal. The analysis of this curve in a Log-Log domain becomes a straigth line.

$$\log C(r) = \log \alpha + D \log r$$
(10.)

The slope of this line is the Fractal dimension D. Since a fractal is self similar, a fractal dimension can be evaluated by comparing a property between any different two scales. In practice D can be estimated by fitting the equation 10. over a range of r.

$$D = 2 - H$$
(11.)

A fractional Brownian motion (fBm) is a continuous-time Gaussian process depending on the Hurst parameter $0 < H < 1$. It generalizes the ordinary Brownian motion corresponding to H = 0.5 and whose derivative is the white noise. The fBm is self-similar in distribution and the variance of the increments is given by

$$\mathrm{Var}(\mathrm{fBm}(t) - \mathrm{fBm}(s)) = v |t - s|^{2H}$$
(12.)

where v is a positive constant.

According to the value of H, the FBm exhibits for H > 0.5, long-range dependence and for H < 0.5, short or intermediate dependence. The Figure 7 shows a piecewise FBm with first of 64 first points H=0.3 and next every 64 points H=H+0.2 until H=0.9. We try to extract the local difference of complexity. H=0.9 clearly exhibits a stronger low-frequency component and has, locally, less irregular behavior than H=0.3.
This signal should exhibit a medium value of the Hurst exponent using classical estimation.
The first method called regularisation method (using Fraclab) estimates D=1.4881, the second computed with fractal dimension of waveforms algorithm [7] found D=1.4892. Therefor we can calculate with Eq 11. the value of H close to H=0.52. It means normal approach only exhibit the global complexity content.

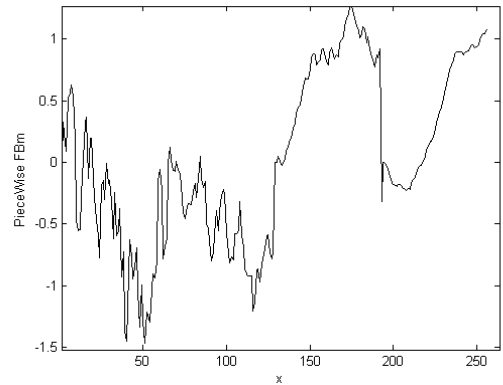

**Figure** 7. Piecewise FBm

This simple example will permits to easily estimate the local behaviour of such FBm and should be compare to wavelet based method [9,10].

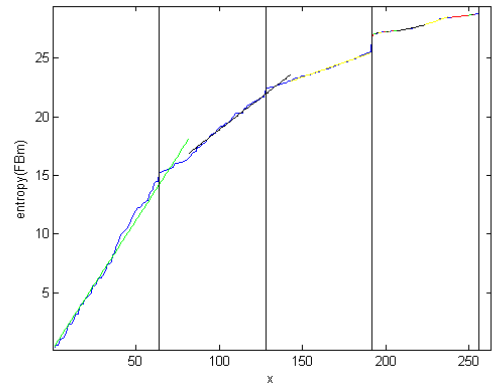

**Figure** 8. Detected lines of the entropy(FBm)

The complexity of the signal can explain the errors of regression. However the method is able to segment the piecewise FBm in 4 parts with sucess. The choice of Rm²= 0.988 is determinant in the result of the procedure.

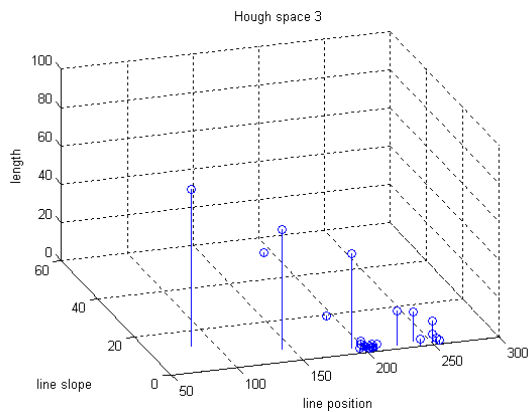

**Figure** 9. Hough parameters space of the FBm

The general model used to fit the variation of tangent function on the Hurst exponent (Figure 10) is f(x) = a*exp(b*x). The estimated coefficients (with 95% confidence bounds) are:
a =   0.6374  (0.5104, 0.7643)
b =   -3.541  (-4.063, -3.018)

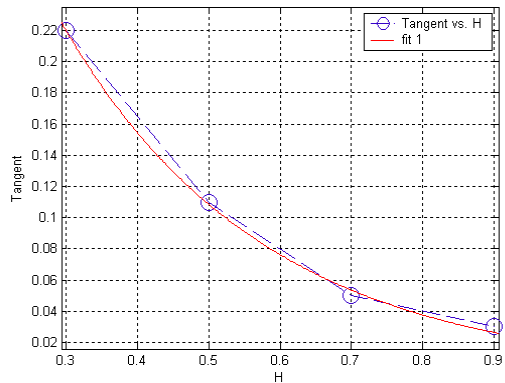

**Figure** 10. Hurst exponent H Vs Tangents of detected lines

The results demonstrates that Hurst exponent Vs Tangents of detected lines can be modeled with an exponential fit with high confidence (Goodness of fit, $R^2$=0.9987). Thus theorical aspects of this method should be studied with attention in order to correlate the Hurst exponent with the slope (or Tangent) of detected lines from the entropy of FBm.

Finally we show as second example a classical mechanical vibration example to verify the performance of the algorithm. This is a simulated modeshape of a damaged cantilever beam (Mode 1 of vibration computed using Finite Element Analysis). The damage is localised at point 20 (on total length 60). Then it becomes difficult to recognize the damage only scaning modeshape data. This is the reason why some important researchs focused on modeshape damage detection using curvature [7,8], wavelets [4,7] and fractal dimension [7]. The entropy function is estimated using discrete formula (6) and plotted with the experimental signal on Figure (4). Choice of $Rm^2$= 0.999.

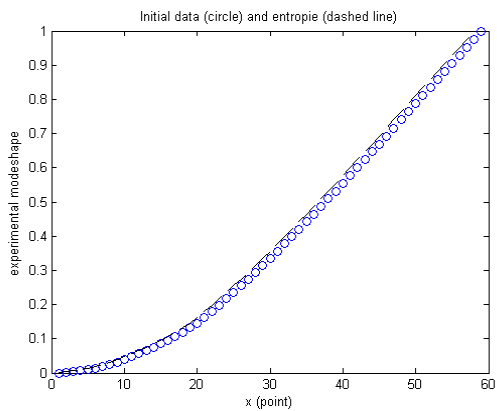

**Figure** 11. Experimental modeshape and entropy function

The detected lines are plotted in the Hough parameters space (Figure 12), it is easy to see there is an abrupt changes in the detected line parameter around point 20 as its length decrease and its slope has a non linear behaviour.

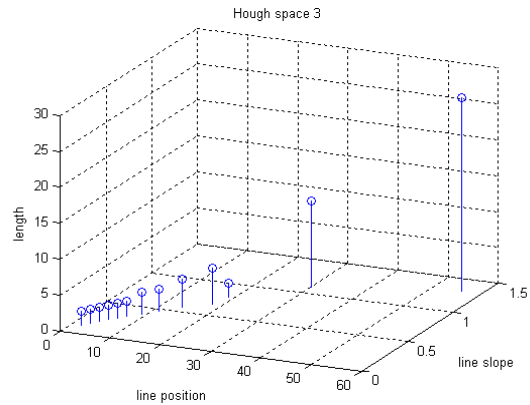

**Figure** 12. The singularity behaviour is represented by a quick change in the Hough parameters space at line position 20

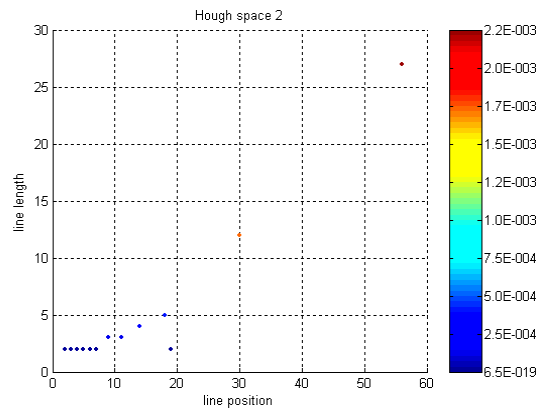

**Figure** 13. Abrupt changes in the detected line length permits to recognize the singularity.

Moreover the Figure 13 enhances the damage with a loss of length in the detected line before point 20. The small lines detected before the damage represents the curvature of the beam (closed to the encastred boundary) with very small error of regression ( error blue in colorbar), the signal becomes more homogeneous after the damage as the length parameter increase (error red in colorbar).

## 6   Conclusions

A new method of detecting homogeneous zones of a 1D signal is proposed. This 1D pattern recognition algorithm has also the abilitie to recognize singularity (or abrupt changes) of the signal. The simple mathematical transformation called entropy permits to study any signal as a piecewise linear one. The interest of the multiple regressions algorithm is its low time computing for extracting tangents and intersepts of detected lines. Plotting these parameters in the Hough parameters space, it becomes simple to interpret the classification results. Finally two experimental validations are proposed, Firstly in order to characterise a piecewise FBm signal and seconly to identify damages on beam using dynamic modeshape data. Further works should be done to automate the choice of R² developing an adaptative multiple regression algorithm. Indeed the choise of the R² (or error between model and fit) is the only user input but require human expert as tool to conclude on the final signal segmentation. The method should be extented to other signal for example exhibit the local frequency shift of a chirp with localised fit of sinusoids.

## 7  Aknowledgment

The author express his thanks to Mr. Cremoux for its comments and also the anonymous referees from the Matlab central website for the imporvement of the paper.

## 8  References

1. A. Denis, F. Cremoux, Using the Entropy of Curves to Segment a Time or Spatial Series, Mathematical Geology, November 2002, vol. 34, no. 8, pp. 899-914(16)
2. G. Ferrari-Trecate and M. Muselli. A new learning method for piecewise linear regression. International Conference on Artificial Neural Networks (ICANN02), 2002. Madrid, Spain.
3. F. Gritzali and G. Papakonstantinou, Fast piecewise linear approximation algorithm, Signal Processing, Volume 5, Issue 3, May 1983, Pages 221-227.
4. J.C. Hong, Y.Y. Kim, H.C. Lee and Y.W. Lee, Damage detection using the Lipschitz exponent estimated by the wavelet transform: applications to vibration modes of a beam, International Journal of Solids and Structures, Volume 39, Issue 7 (2002), p. 1803-1816.
5. H.E. Hurst, Long-term storage capacity of reservoirs. Trans.Am.Soc. Civil Eng. 116: 770-808 (1951).
6. Atsushi Imiya, Detection of piecewise-linear signals by the randomized Hough transform, Pattern Recognition Letters,Volume 17, Issue 7 , 10 June 1996, Pages 771-776.
7. J. Morlier , F. Bos , P. Castera, Benchmark of damage localisation algorithms using mode shape data, Proceedings of the 6th International Conference on Damage Assessment of Structures, july 2005.
8. H.B. Nugraha and A.Z.R. Langi, "Segmented fractal dimension measurement of 1-D Signal: A wavelet based method", 2002 Asia Pacific Conference on Communications (APCC 2002), Bandung, Indonesia, 17-19 September 2002.
9. H.B. Nugraha and A.Z.R. Langi, "A Wavelet-Based Measurement of Fractal Dimensions of a 1-D Signal", ", International Conference on Information, Computer, and Signal Processing (ICICS 2001), Mandarin Hotel Singapore, 15-18 Oct 2001.
10. A.K. Pandey, M. Biswas and M.M. Samman: Damage Detection from Changes in Curvature Mode Shapes, J. of Sound and Vibration Vol. 145(2) (1991), p. 321-332.
11. Yujing, Zeng, Janusz Starzyk, Piecewise Linear Approach: a New Approach in Automatic Target Recognition,  Proceedings of SPIE, 2000.
12. T. Risse, Hough Transform for Line Recognition, Computer Vision and Image Processing, 1989, 46, 327-345, 1989.

.